\documentclass[letterpaper]{article} % DO NOT CHANGE THIS
\usepackage{aaai23}  % DO NOT CHANGE THIS
\usepackage{times}  % DO NOT CHANGE THIS
\usepackage{helvet}  % DO NOT CHANGE THIS
\usepackage{courier}  % DO NOT CHANGE THIS
\usepackage[hyphens]{url}  % DO NOT CHANGE THIS
\usepackage{graphicx} % DO NOT CHANGE THIS
\usepackage{natbib}  % DO NOT CHANGE THIS AND DO NOT ADD ANY OPTIONS TO IT
\usepackage{caption} % DO NOT CHANGE THIS AND DO NOT ADD ANY OPTIONS TO IT
\usepackage{algorithm}
\usepackage{algorithmic}
\usepackage{amsmath}
\usepackage{amssymb}
\usepackage{mathtools}
\usepackage{amsthm}
\usepackage{ulem}
\usepackage{amsfonts}
\usepackage{comment}
\usepackage{epstopdf}
\usepackage{algorithm}
\usepackage{mathtools}
\mathtoolsset{showonlyrefs=true}
\usepackage{color}
\usepackage{hyperref}
\usepackage{ulem}
\usepackage{natbib}

\usepackage{pgfplots}
\pgfplotsset{compat = newest}
\usepgfplotslibrary{colormaps}

\usepackage[capitalize,noabbrev]{cleveref}

\newcommand{\red}[1]{\textcolor{red}{#1}}

\theoremstyle{plain}
\newtheorem{theorem}{Theorem}[section]

\newtheorem{lemma}[theorem]{Lemma}
\newtheorem{corollary}[theorem]{Corollary}
\theoremstyle{definition}
\newtheorem{definition}[theorem]{Definition}
\newtheorem{assumption}[theorem]{Assumption}
\newtheorem{remark}[theorem]{Remark}
\newtheorem{example}[theorem]{Example}
\newtheorem{question}[theorem]{Question}

\usepackage{newfloat}
\usepackage{listings}
\DeclareCaptionStyle{ruled}{labelfont=normalfont,labelsep=colon,strut=off} % DO NOT CHANGE THIS
\floatstyle{ruled}
\newfloat{listing}{tb}{lst}{}
\floatname{listing}{Listing}
\title{Linear Regularizers Enforce the Strict Saddle Property}
\title{Linear Regularizers Enforce the Strict Saddle Property}
\author {
    Matthew Ubl,\textsuperscript{\rm 1}
    Kasra Yazdani, \textsuperscript{\rm 1}
    Matthew T. Hale \textsuperscript{\rm 1}
}
\affiliations {
    \textsuperscript{\rm 1} Department of Mechanical and Aerospace
Engineering \\
University of Florida, Gainesville, FL, 32611, USA. \\
    m.ubl@ufl.edu, kasra.yazdani@ufl.edu, matthewhale@ufl.edu
}
\begin{document}

\maketitle

\begin{abstract}
Satisfaction of the strict saddle property has become a standard assumption in non-convex optimization, and it ensures that many first-order optimization algorithms will almost always escape saddle points. However, functions exist in machine learning that do not satisfy this property, such as the loss function of a neural network with at least two hidden layers. First-order methods such as gradient descent may converge to non-strict saddle points of such functions, and there do not currently exist any first-order methods that reliably escape non-strict saddle points. To address this need, we demonstrate that regularizing a function with a linear term enforces the strict saddle property, and we provide justification for only regularizing locally, i.e., when the norm of the gradient falls below a certain threshold. We analyze bifurcations that may result from this form of regularization, and then we provide a selection rule for regularizers that depends only on the gradient of an objective function. This rule is shown to guarantee that gradient descent will escape the neighborhoods around a broad class of non-strict saddle points, and this behavior is demonstrated on numerical examples of non-strict saddle points common in the optimization literature.
\end{abstract}

\section{Introduction}
\label{intro}
Interest in non-convex optimization has grown in recent years, driven by applications such as 
training deep neural networks. 
Often, one seeks convergence to a local minimizer in such problems because finding global minima is known to
be NP complete~\cite{murty87}. To ensure convergence to minimizers, one research direction 
in non-convex optimization has been the identification of problem properties for which particular algorithms escape saddle points. One such property, which has become common in the non-convex optimization literature since its introduction in~\cite{ge2015escaping}, is the \textit{strict saddle property} (SSP), which states that the Hessian of every saddle point of a function has at least one negative eigenvalue. It was later shown that gradient descent and other first order methods almost always escape saddle points of objective functions
that satisfy the SSP (and other mild assumptions)~\cite{lee2016gradient,panageas2017gradient,lee2017first}. 

Because of this behavior, a growing body of non-convex optimization research has 
either focused on problems for which the SSP is known to hold, or simply assumed the SSP holds for a 
generic problem and derived convergence guarantees that result from it. 
However, verification of the SSP for a general, unstructured problem is difficult in practice, 
and there exist problems in machine learning for which the SSP does not hold, such as training a neural 
network with at least two hidden layers~\cite{kawaguchi2016deep}. 

Motivated by these challenges, we develop 
%an algorithmic method that will
a linear regularization framework that will allow first-order methods to escape saddle points that are not 
strict. Specifically, our approach is to \textit{enforce} the SSP
by regularizing problems when in the vicinity of a non-strict saddle point, rather than simply assuming that the SSP holds. 
We show that this can be done with a linear regularizer, 
motivated by John Milnor's proof that almost all choices of such a term will 
render a function Morse (and therefore 
enforce the SSP)~\cite{milnor2015lectures}.
We are also motivated by the success of regularization techniques in convex optimization, where quadratic
perturbations are used to provide strong convexity to objective functions~\cite{facchinei07}, and we believe that
the linear regularizers we present are their natural counterparts in the non-convex setting. 

\subsection{Related Work}
\label{intro.relatedwork}

A large body of work exists on the convergence properties of gradient descent and other first-order methods on problems with the SSP, including algorithms that 
consider deterministic gradient descent~\cite{dixit2020exit,schaeffer2019extending}, and those that incorporate noise into their updates~\cite{xu2017first,daneshmand2018escaping,yang2017fast,ge2015escaping}. These methods are shown to escape strict saddles, but have not been shown to escape non-strict saddles, and therefore rely on the SSP.

%incorporate noise into their updates~\cite{xu2017first,daneshmand2018escaping,yang2017fast,ge2015escaping} and those that consider deterministic gradient descent~\cite{dixit2020exit,schaeffer2019extending}. Zeroth order methods have also been developed~\cite{bai2020escaping}. All of these methods guarantee convergence to a second-order optimal point, i.e., a critical point with a positive semi-definite Hessian, and therefore need the SSP to ensure such a point is a local minimum (otherwise they may be non-strict saddles).

While these methods are shown to escape strict saddles in the limit, they can get stuck near strict saddles for exponential time, which can cause numerical slowdowns~\cite{du2017gradient}. Attempts have been made to accelerate the escape near strict saddle points~\cite{jin2017escape,agarwal2017finding,jin2018accelerated}. However, first-order methods may actually converge to non-strict saddles, and such accelerated methods do not escape.

Current research into escaping non-strict saddle points uses higher-order information and/or algorithms. 
Perhaps the best known is~\cite{anandkumar2016efficient}, which guarantees convergence to a third-order optimal critical point. That paper replaces the SSP, which is a property of the Hessian, with a condition on the third-order derivative of the objective function. Work in~\cite{zhu2020adaptive} expands on these results and includes simulations for a function that does not satisfy the SSP. Later work in~\cite{chen2021high} provides a method to converge to $p^{th}$-order critical points using $p^{th}$-order information, while also demonstrating that doing so is NP-hard for $p \geq 4$. Recent work in~\cite{truong2021new} examines the behavior of a second-order method on common examples of non-strict saddle points, and~\cite{nguyen2017loss} develop a weaker form of the SSP that guarantees escape from saddle points when training a particular neural network. In contrast, we require only first-order information
and provably escape from non-strict saddles using linear regularizers under
weak assumptions. 

Previous research has shown that regularizing with quadratic or sums of squares (SOS) terms will make a function Morse, which is sufficient to ensure the SSP is satisfied~\cite{lerario2011plenty,nicolaescu2011invitation}. However, no convergence or bifurcation analysis was performed on the regularized function, and
indeed these 
results originate outside the non-convex optimization literature. 
We show in Example~\ref{e.quadregbad} that quadratic and SOS
regularizers can actually convert a non-strict saddle point into a local minimum, 
and thus we do not use them.

\subsection{Contributions}
\label{intro.contributions}

The contributions of this paper are the following: 
\begin{itemize}
    \item We identify certain properties that any linear regularization scheme must have, namely that regularizers cannot be chosen randomly, must be chosen locally, and must have their norms obey an upper bound dependent on $f$.
    \item We present a regularization scheme that has the above properties, and analyze the bifurcations it induces.
    \item We prove that, under a condition much weaker than the SSP, the presented regularization scheme escapes all saddle points (strict and non-strict) of $f$.
    \item We bound the regularization error seen at minima that is induced by linear regularizers.
\end{itemize}

The remainder of the paper is organized as follows. Section~\ref{linearregularizers} establishes the theoretical motivation behind a linear regularization scheme. 
In Section~\ref{bifurcations}, we analyze the bifurcations that may occur when regularizing, 
identify the properties a linear regularization scheme for SSP enforcement must have, and present a particular choice of regularizer that has these properties. 
In Section~\ref{exitcondition}, we prove this regularization 
method escapes saddle points that satisfy a condition weaker than the SSP and demonstrate this escape on examples of non-strict saddle points taken from the literature. 
In Section~\ref{hyperparameter}, we analyze a hyperparameter that regulates the size of regularization and its effect on speed and accuracy, and in Section~\ref{concludingremarks} we provide concluding remarks.

\section{Linear Regularization} \label{linearregularizers}

Throughout this paper, $f: \mathbb{R}^n \rightarrow \mathbb{R}$ denotes a function in $C^2$, the space of twice-continuously differentiable functions, with $L$-Lipschitz gradient $\nabla f$. 
The symbol $g: \mathbb{R}^n \rightarrow \mathbb{R}^n$ denotes
a first-order map, with iterates generated by the sequence $x_k = g(x_{k-1}) = g^k(x_0)$. For clarity, in this paper we take $g$ to represent a gradient descent mapping, i.e., $g(x) = x - \gamma \nabla f(x)$, with $\gamma \in (0, 1/L)$, though we note the results of this paper hold for any choice of $g$ that avoids strict saddle points, see~\cite{lee2017first}. 
The following definition regards the \textit{critical points} of $f$:

\begin{definition} 
\begin{enumerate} 
    \item A point $x^*$ is a \textit{critical point} of $f$ if $\nabla f(x^*) = 0$ or, equivalently, $g(x^*) = x^*$. 
    \item A critical point $x^*$ is \textit{isolated} if there exists a neighborhood $U$ around $x^*$ with $x^*$ as the only critical point in $U$. Otherwise it is called \textit{non-isolated}. 
    \item A critical point of $f$ is a \textit{local minimum} if there exists a neighborhood $U$ around $x^*$ such that $f(x^*) \leq f(x)$ for all $x \in U$, and a \textit{local maximum} if $f(x^*) \geq f(x)$. 
    \item A critical point of $f$ is a \textit{saddle point} if for all neighborhoods $U$ around $x^*$, there exist $y,z \in U$ such that $f(y) \leq f(x^*) \leq f(z)$. 
    \item A critical point of $f$ is a \textit{strict saddle} if $\lambda_{min}(\nabla^2f(x^*)) < 0$. 
    \item The \textit{local stable set} $W^{s}_{g} (x^*)$ defined on some neighborhood $U$ 
    of a critical point $x^*$ is the set of initial conditions of the first-order map $g$ in $U$ that converge to $x^*$, i.e., 
    $W^s_{g} (x^*) = \{x \in U:\lim_{k \to \infty} g^k(x)=x^*\}$. 
    The \textit{local unstable set} is defined as $W^u_{g} (x^*) = \{x \in U:\lim_{k \to \infty} g^k(x)\neq x^*\}$. 
    If $U = \mathbb{R}^n$, then $W^s_{g} (x^*)$ ($W^u_{g} (x^*)$) is the \textit{global stable (unstable) set}. 
\end{enumerate}
\end{definition}

Here~$\lambda_{min}(\cdot)$ denotes the minimum eigenvalue of a square matrix. 
Lemma~\ref{l.gdavoidsstrictsaddles} states that, for almost all initial conditions, $g^k(x)$ does
not converge to a strict saddle:

\begin{lemma} \label{l.gdavoidsstrictsaddles}
\cite{panageas2017gradient}
Let $f : \mathbb{R}^n \rightarrow \mathbb{R}$ be a $C^2$ function with $L$-Lipschitz gradient. The set of initial conditions $x \in \mathbb{R}^n$ such that $g^k(x)$ converges to a strict saddle point of $f$ is of (Lebesgue) measure zero.
\end{lemma}
\textit{Proof:} See Theorem 2 in~\cite{panageas2017gradient}. $\hfill\square$

%The intuition behind this result comes from the fact that for any strict saddle $x^*$ of $f$, $W^s(x^*)$ is measure zero on $\mathbb{R}^n$. 
%This is visualized in Example~\ref{e.strictsaddle}: 
%
%\begin{example} \label{e.strictsaddle}
%\textit{(Measure Zero $W^s$ for Strict Saddle)}
%
%Consider the function $f(x,y) = \frac{1}{2}x^2 - \frac{1}{2}y^2$, with streamlines plotted on Figure~\ref{fig:strictsaddle}. Here, $(0,0)$ is a strict saddle of $f$, with $\nabla^2 f(0,0)$ having $1$ and $-1$ as eigenvalues. The positive eigenvalue corresponds to positive curvature along the $x$-axis, and the negative eigenvalue corresponds to
%negative curvature along the $y$-axis. $W^s(0,0)$ is represented by the red area and $W^u(0,0)$ by green. It is easy to see that gradient descent only converges to $(0,0)$ from initial conditions on the~$x$-axis, 
%because any others will escape along the direction of negative curvature along the $y$-axis.
%\end{example}
%
%\begin{figure}[htp]
%    \centering
%    \includegraphics[width=6cm]{StrictSaddleSketch.jpeg}
%    \caption{A streamline plot for gradient descent applied to the function~$f(x, y) = x^2-y^2$. 
%    The pink region is all initial conditions that converge to $(0,0)$ and the green region is all
%    initial conditions that diverge.}
%    \label{fig:strictsaddle}
%\end{figure}

%While the exact topology of strict saddles of other functions may be more complex, 
The underlying principle is that,
for a saddle~$x^*$,
a single negative eigenvalue 
of $\nabla^2 f(x^*)$ renders $W^s_{g}(x^*)$ measure zero. This is the motivating 
principle behind the study of the strict saddle property:

\begin{definition}
A function $f$ satisfies the \textit{strict saddle property (SSP)} if every saddle point of $f$ is strict.
\end{definition}

From Lemma~\ref{l.gdavoidsstrictsaddles}, gradient descent will almost always avoid every strict saddle point of an objective function $f$.
Therefore, if $f$ satisfies the SSP, then gradient descent will almost always avoid \textit{all} saddle points of $f$. Provided gradient descent converges (i.e., $\lim_{k \to \infty} g^k(x)$ exists), it must then 
almost always 
converge to a local minimum. We note that $g^k(x)$ is guaranteed to converge in a variety of settings, including when $f$ is analytic or coercive, %such as when $f$ is a loss function, 
and we will proceed with the assumption that $f$ satisfies one of these properties.

However, verifying that a general, unstructured function satisfies the 
SSP is difficult in practice, and functions of interest exist that 
are known not to satisfy the SSP, such as the loss function of a neural network with at least two hidden layers~\cite{kawaguchi2016deep}. These functions may have \textit{non-strict saddles}:

\begin{definition}
A saddle point $x^*$ of $f$ is a \textit{non-strict saddle} if $\lambda_{min}(\nabla^2f(x^*)) = 0$.
\end{definition}

We make a brief point on terminology here. The definition of a degenerate saddle 
varies between the dynamical systems and computer science literature, so to avoid 
confusion in this paper a \textit{degenerate saddle} is any saddle point $x^*$ 
whose Hessian has at least one zero eigenvalue (i.e., $\nabla^2 f (x^*)$ is singular), 
while a \textit{non-strict saddle} is a saddle with a  Hessian whose minimum eigenvalue is zero (i.e., $\nabla^2 f (x^*)$ is singular \textit{and} positive semi-definite). Using this terminology, any non-strict saddle is necessarily degenerate. We note that 
the SSP is \textit{not} a non-degeneracy condition, as the Hessians of strict saddles may be degenerate, as long as they have at least one negative eigenvalue. Example~\ref{e.degeneratesaddle} illustrates 
the key problem with non-strict saddle points, which is 
that their stable sets are not necessarily measure zero.

\begin{figure}[htp]
    \centering
    \includegraphics[width=4cm]{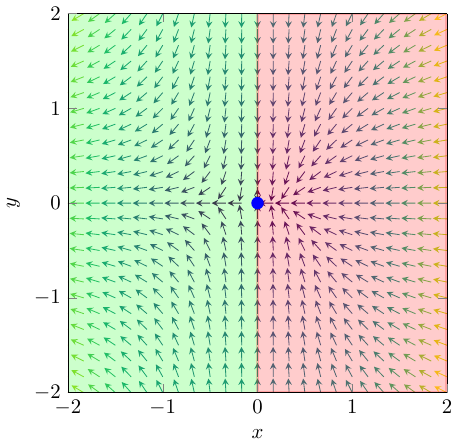}
    \caption{The negative gradient field of $f(x,y) = \frac{1}{3}x^3+\frac{1}{2}y^2$. The blue dot at $(0,0)$ denotes the non-strict saddle point, $W^u_{g}(0,0)$ is denoted by the green region, and $W^s_{g}(0,0)$ by red.}
    \label{fig:degeneratesaddle}
\end{figure}

\begin{example} \label{e.degeneratesaddle}
%\textit{(Non-Measure Zero $W^s_{\notes{g}}$ for non-Strict Saddle)} 
Consider the function $f(x,y) = \frac{1}{3}x^3 + \frac{1}{2}y^2$, with negative gradient field plotted in Figure~\ref{fig:degeneratesaddle}. Here, $(0,0)$ is a non-strict saddle of $f$, with $\nabla^2 f(0,0)$ having $1$ and $0$ as eigenvalues. We see that $W^s_{g}(0,0) = \{(x,y) : x > 0\}$, depicted by the red region. That is, the set of initial conditions for which $g^k(x,y)$ converges to $(0,0)$ is not measure zero and is in fact a closed halfspace of $\mathbb{R}^2$.
\end{example}

Instead of modifying gradient descent to somehow accommodate non-strict saddles, we instead wish to modify the \textit{problem} itself in such a way that the modified function satisfies the SSP, either by making non-strict saddles strict or eliminating them altogether. That is, we wish to find a regularization scheme that enforces satisfaction of the SSP and thus
ensures the escape of non-strict saddles. While quadratic and sums of squares regularizers are used in convex optimization, 
they can be harmful in non-convex problems
because they can
change the positive semi-definite Hessian of a non-strict saddle into a positive definite one, turning such a saddle into a local minimum: 

\begin{example} \label{e.quadregbad}
Consider again the function $f(x,y) = \frac{1}{3}x^3 + \frac{1}{2}y^2$, which has a non-strict saddle at $(0,0)$ with eigenvalues $1$ and $0$. If a sum of squares regularization term $\frac{1}{2}\alpha_x x^2 + \frac{1}{2}\alpha_y y^2$  is added to $f$, then $(0,0)$ remains a critical point of the regularized function, but the eigenvalues of the regularized Hessian become $\alpha_x$ and $1+\alpha_y$, rendering $(0,0)$ a local minimum for all $\alpha_x,\alpha_y > 0$.
\end{example}

Instead, the following lemma provides motivation for using a linear regularization term.

\begin{lemma} \label{l.milnor}
\cite{milnor2015lectures}
If $f : \mathbb{R}^n \rightarrow \mathbb{R}$ is a $C^2$ function, then for almost all $l \in \mathbb{R}^n$, 
the critical points of the function $f_l(x) = f(x) + l^T x$ have only non-singular Hessians.
\end{lemma}
\textit{Proof:} See Lemma~A in~\cite{milnor2015lectures}. $\hfill\square$

This lemma states that for almost any choice of $l$ (any except a set of Lebesgue measure zero) the regularized function $f_l$ will have only non-degenerate critical points. The fact that non-degenerate saddles are strict immediately gives us the following corollary:

\begin{corollary} \label{t.enforcessp}
If $f : \mathbb{R}^n \rightarrow \mathbb{R}$ is $C^2$, then for almost all $l \in \mathbb{R}^n$, the function 
$f_l(x) = f(x) + l^T x$ satisfies the SSP.
\end{corollary}
%\textit{Proof:} This follows from Lemma~\ref{l.milnor} and the fact that non-degenerate saddle points must be strict. $\hfill\square$

This regularization method does not affect the Hessian (i.e., $\nabla^2 f(x) = \nabla^2 f_l(x)$), avoiding the problems
caused by sums of squares and quadratic regularizers. 
Corollary~\ref{t.enforcessp} now motivates the following question, which will be the focus of the remainder of this paper:
\begin{question} \label{q}
Can a linear regularization scheme be used to enforce the SSP on functions that do not satisfy it? If so, what properties must such a scheme have?
\end{question}
%This suggests that linear regularization of a non-convex problem may be a natural counterpart to quadratic regularization of convex problems.
%\notes{Corollary~\ref{t.enforcessp} answers an open question from~\cite{lee2016gradient} 
%that asks ``\textit{Will a perturbation of a
%function always satisfy the strict saddle property?}'' in the affirmative.} 
Though Corollary~\ref{t.enforcessp} states that almost every choice of $l$ will enforce the SSP, it is important to understand \textit{how} the SSP is enforced. As we will see 
in the following section, this regularization method enforces satisfaction of the SSP by creating bifurcations of degenerate critical points of $f$, and we must carefully analyze these bifurcations to ensure
that we attain the desired convergence properties.

\section{Bifurcations}
\label{bifurcations}

\begin{figure}[htp]
    \centering
   \includegraphics[width=8cm]{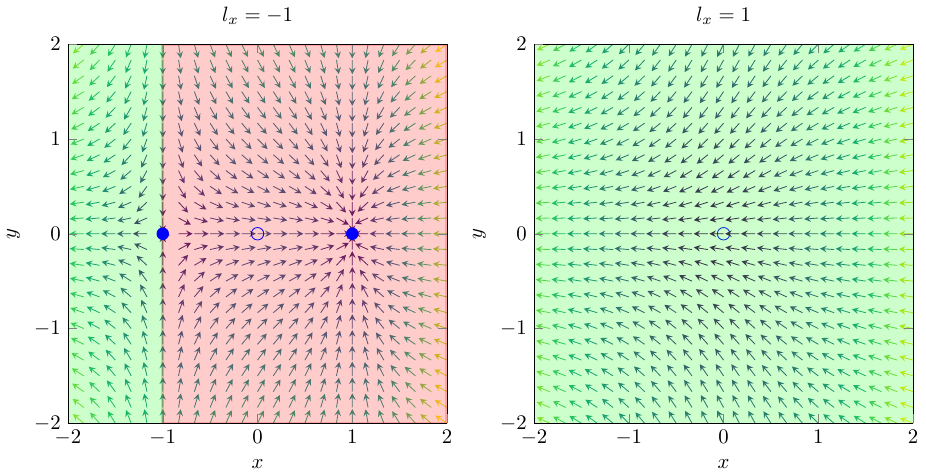}
    \caption{With $l_x = -1$ we create a local minimum \textit{and} a strict saddle, and the escape region shifts (left). 
    With $l_x = 1$, the critical point is destroyed and gradient descent escapes the saddle from every initial condition (right).}
    \label{fig:degeneratesaddlebifurcation}
\end{figure}

Regularization of a function perturbs non-degenerate critical points, which can be limited by a judicious choice
of regularizer. 
However, the same is not true of degenerate critical points, as can be seen in the following example.
%However, more seriously, regularization can also cause a \textit{local bifurcation}  at a degenerate critical point:
\begin{example} \label{e.bifurcation}
Consider again the function ~$f(x, y) = \frac{1}{3}x^3 + \frac{1}{2}y^2$ and consider two regularizations that add terms 
of the form $l_x x + l_y y$. The first sets~$l_x = 1$ and~$l_y = 0$ and
the second sets~$l_x = -1$ and~$l_y = 0$, and we
plot the trajectory behavior of gradient descent for each in Figure~\ref{fig:degeneratesaddlebifurcation}.
\end{example}

Observe that when $l_x = -1$, the original non-strict saddle splits 
into a strict saddle at $(-1,0)$ and a local minimum at $(1,0)$. Both of these 
points are non-degenerate, satisfying the SSP as ensured by Corollary~\ref{t.enforcessp}. 
However, we can see that $W^s_g$
(now defined for both of the resulting critical points, shown in red) has actually expanded. We have observed a \textit{local bifurcation} of the non-strict saddle point at $x^*$.

\begin{definition}
Let $h: \mathbb{R}^n \times \mathbb{R}^k \rightarrow \mathbb{R}$ be a $C^2$ function. Let $(x^*,\mu^*)$ be a point for which $\nabla_x h(x^*,\mu^*) = 0$ and $\nabla_x^2 h(x^*,\mu^*)$ is singular. A \textit{local bifurcation} of this gradient system occurs at $x^*$ when a smooth change in the parameter $\mu$ away from $\mu^*$ induces a sudden change in the stability properties of the negative gradient vector field at $x^*$.
%\mhmargin{Is there a concise way to say what object's stability properties are the ones changing here?}
\end{definition}

A ``sudden change in stability properties'' can mean
a number of things, see~\cite{guckenheimer2013nonlinear}, but in the situation presented in this paper (a codimension-one linear perturbation of a gradient system) it refers almost exclusively to \textit{saddle-node bifurcations}. 
Example~\ref{e.bifurcation}, for which $h(x,y,\mu) = \frac{1}{3}x^3+\frac{1}{2}y^2+\mu x$, illustrates
a saddle-node bifurcation, where a degenerate critical point at $x^*$ splits into two or more critical points, or the critical point at $x^*$ is eliminated. 
This bifurcation occurs when~$\mu$ crosses from zero to being positive or negative, and it
results in $W^s_g(x^*)$ changing size or dimension. Note that the saddle-node bifurcation in Example~\ref{e.bifurcation} has created a \textit{false minimum} at $(1,0)$:

\begin{definition} \label{d.falsemin}
A \textit{false minimum} is a local minimum of $f_l$ that resulted from a bifurcation of a degenerate saddle point of $f$
that was caused by the linear regularizer~$l^Tx$. 
\end{definition}

In Example~\ref{e.bifurcation}, one can see that for any $l_x < 0$, a saddle-node bifurcation occurs. 
We also observe that when $l_x = 1$ (and in fact whenever $l_x > 0$) the critical point 
at $(0,0)$ is destroyed and all trajectories 
of gradient descent 
escape the neighborhood of $(0,0)$ (i.e., $W^u_g = \mathbb{R}^2$, shown in green). This gives us the following remark regarding Question~\ref{q}: 
\begin{remark} \label{r.notrandom}
Any linear regularization scheme that chooses $l$ randomly has a positive probability of creating a false minimum near a non-strict saddle point of $f$.
%\mhmargin{This is true, but can you make it easier on the reader? The ``random'' part comes out of nowhere. This could likely be fixed by adding
%a sentence before: something like ``If you just draw one out of thin air, you might make one of these bad things (that you just saw) happen, so now we have this result''.
%}
\end{remark}
Intuitively then, $l$ should have some dependence on~$f$, and 
$\nabla f$ specifically is the only information available to a first-order algorithm. We note that because $l$ cannot be chosen randomly, we cannot rely solely on Corollary~\ref{t.enforcessp} to guarantee that a particular choice of $l$ enforces the SSP. 
%\red{Determining a selection rule for $l$ and proving it satisfies the SSP will be the focus of the reminder of this paper.} 
%\mhmargin{There was a sentence like this above Question 2.9. It's just jarring for the reader to see this twice so close together. I'm fine with cutting this one if you are -- the one above Question 2.9 is really useful so that one should stay.} 

We present the following example to illustrate another property a linear regularization scheme must have.

%This example demonstrates 
%that while almost all choices of $l$ will enforce the SSP, only some enforce the trajectory behavior 
%we want. Thus, $l$ should be chosen deliberately, based on some information from the function being minimized, 
%rather than arbitrarily. However, a choice of $l$ that enforces desired behavior around one non-strict saddle 
%point may not do so around another, as the following theorem states.

\begin{example} \label{t.timevaryingbad}
The function $f(x) = (x-1)^3(x+1)^3$ has non-strict saddles at $x = -1$ and $x = 1$. For any arbitrarily small choice of $l > 0$, the non-strict saddle at $x = -1$ undergoes a saddle-node bifurcation and the non-strict saddle at $x = 1$ is destroyed. For any arbitrarily small choice of $l < 0$, the non-strict saddle at $x = 1$ experiences a saddle-node bifurcation and the non-strict saddle at $x = -1$ is destroyed.
\end{example}

\begin{figure}[htp]
    \centering
    \includegraphics[width=8cm]{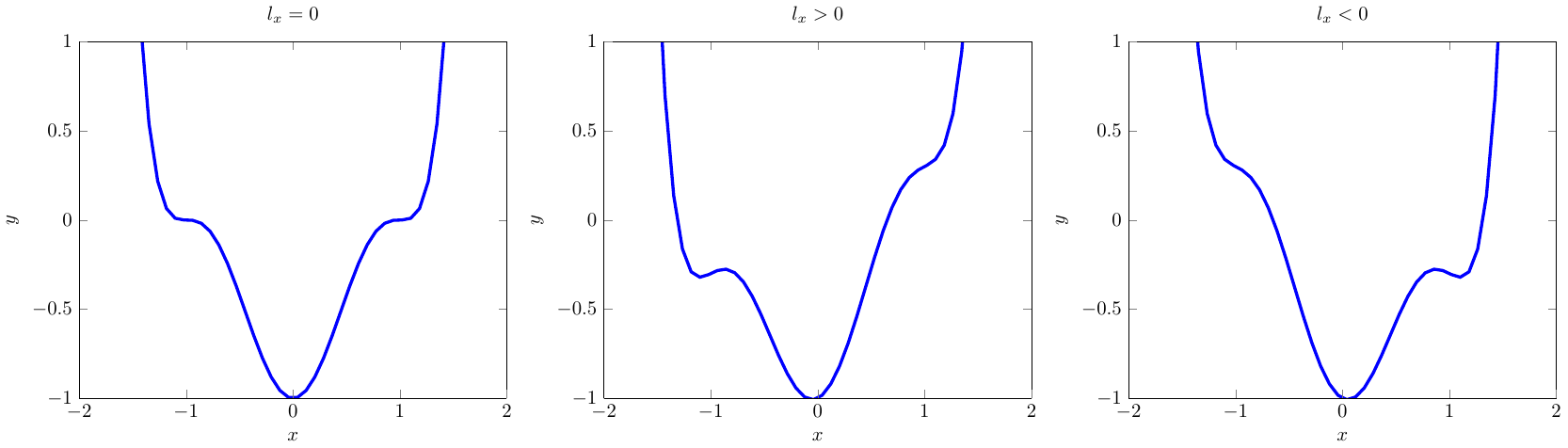}
    \caption{
    Plots of the function
    $(x-1)^3(x+1)^3 + lx$ for~$l = 0$, $l > 0$, and~$l < 0$. 
    Regardless of the sign of $l$, one of the original degenerate critical points is bifurcated into a false minimum and a local maximum, and the other is eliminated
    for every regularizer~$l \neq 0$.
}
    \label{fig:timevarying}
\end{figure}

A natural consequence of Example~\ref{t.timevaryingbad} is the following remark regarding Question~\ref{q}:
\begin{remark} \label{r.notglobal}
There exist $C^2$ functions for which any constant, global choice of $l \neq 0$ creates a false minimum.
\end{remark}
Therefore, a linear regularization scheme should choose $l$ ``locally'', changing the choice of $l$ when in the neighborhood of 
different critical points. In order to do so practically we take inspiration from~\cite{jin2017escape} and define a ``small gradient region'', outside of which $l = 0$ and inside of which $l$ will be chosen according to some selection rule that we will devise below:

\begin{definition}
Fix~$\theta > 0$ and
let $L_{\theta}= \{x \in \mathbb{R}^n : \|\nabla f(x) \|_{2} \leq \theta\}$. That is, the \textit{small-gradient region} 
$L_{\theta}$ is the subset of $\mathbb{R}^n$ for which the norm of the gradient of $f$ is less than or equal to $\theta$. For a particular $x \in L_{\theta}$, let the \textit{small-gradient neighborhood} $\Theta (x)$ be the largest connected subset of $L_{\theta}$ that contains $x$.
\end{definition}

As long as $\theta$ is chosen small enough, a point in $L_{\theta}$ must be ``near'' a critical point of $f$.
Local linear regularization means that if an algorithm enters $L_{\theta}$ at some point $x_0$, then the algorithm will choose a regularizer $l$ and use it until it exits $\Theta(x_0)$ (after which~$l$ is reset to zero).  Recall from 
Example~\ref{t.timevaryingbad} that a choice of $l$ that destroys one degenerate critical 
point may induce a saddle-node bifurcation at another. Therefore, 
to avoid a saddle node bifurcation within~$\Theta(x_0)$, 
we must ensure $\Theta(x_0)$ contains at most one critical 
point or connected manifold of critical points. We formalize this idea with the following definition and assumption:

\begin{definition}
Let $X^* = \{x^* \in \mathbb{R}^n : \nabla f(x^*) = 0\}$. That is, $X^*$ is the set of all isolated or non-isolated critical points of $f$. For a particular $x^* \in X^*$, let $\Phi (x^*)$ be the largest connected subset of $X^*$ such that $x^* \in \Phi (x^*)$.
\end{definition}
If $x^*$ is an isolated critical point, then $\Phi (x^*) = \{x^*\}$. If $x^*$ is non-isolated, then $\Phi (x^*)$ is the connected critical manifold that contains $x^*$.

\begin{assumption} \label{a.thetaseparate}
For $f$, there exists $\bar{\theta} > 0$ such that if ${\theta < \bar{\theta}}$, then for every $x^* \in X^*$, $\Theta (x^*) \cap X^* = \Phi (x^*)$.
\end{assumption}

Note that, trivially, $X^* \subset L_{\theta}$ for any $\theta > 0$.
Assumption~\ref{a.thetaseparate} simply states that $\theta$ can be chosen small enough that any critical point $x^*$ is isolated in $\Theta(x^*)$ from all other critical points it is not connected to. 

Recall again from Example~\ref{t.timevaryingbad} that a choice of $l$ that does not induce a saddle-node bifurcation at $x^*$ may do so for other degenerate critical points of $f$. 
We want to ensure that false minima, or indeed any critical points that result from a bifurcation or perturbation of a critical point other than $x^*$, do not end up in the set $\Theta(x^*)$. 
This is guaranteed by the following theorem:

\begin{theorem} \label{t.isolatedbifurcations}
Let $x^*$ be a critical point of $f$, and let $\|l\|_{2} < \theta  < \bar{\theta}$. Let $x^*_l$ be a critical point of the regularized function $f_l$ that resulted as a bifurcation or a perturbation of $x^*$. Then $x^*_l \in \Theta(x^*)$.
\end{theorem}
\textit{Proof:} See Appendix~\ref{apx.t.isolatedbifurcations}. $\hfill\square$

Theorem~\ref{t.isolatedbifurcations} ensures that, even if a particular 
choice of $l$ induces a bifurcation at another degenerate critical point $y^* \in X^*$, the critical 
points that result from that bifurcation are contained within $\Theta(y^*)$, which is disjoint from $\Theta(x^*)$, provided $l$ is sufficiently small. 
%Thus, while the iterates of gradient descent are in~$\Theta(x^*)$, regularizing can create
%false minima somewhere, but not in~$\Theta(x^*)$, and there is no risk of converging to false minima. 
In fact, Theorem~\ref{t.isolatedbifurcations} 
implies that the topology of $\Theta(x^*)$ after regularization depends \textit{only} on the topology 
of $\Theta(x^*)$ prior to regularization. Given this fact, we now wish to choose $l$ such that, if the critical point $x^*$ is a degenerate saddle, 
regularization does not create any false minima in $\Theta(x^*)$. We know from Remark~\ref{r.notrandom} that the choice of $l$ for $\Theta(x^*)$ must depend on the values of $\nabla f$ on $\Theta(x^*)$, and from Theorem~\ref{t.isolatedbifurcations} that we must have $\|l\|_2 \leq \theta$. Upon entering $\Theta(x^*)$ at a point $x_0$, the only value
of~$\nabla f$ over~$\Theta(x^*)$ available is $\nabla f(x_0)$. Therefore it is natural that the choice of $l$ for $\Theta(x^*)$ should be some function of $\nabla f(x_0)$. Two immediate candidates are $l = \nabla f(x_0)$, or $l = -\nabla f(x_0)$. To understand the implications of either of these potential choices, we look at the following theorem:

\begin{theorem} \label{t.3.4.1}
\cite{guckenheimer2013nonlinear} Consider the function $f(x) + \mu l^Tx$ with $\mu \in \mathbb{R}$ and $l,x \in \mathbb{R}^n$. Assume that for $\mu = 0$ there exists a critical point $x^*$ such that:
\begin{enumerate}
    \item $\nabla^2 f(x^*)$ has $n-1$ positive eigenvalues, and a simple eigenvalue 0 with eigenvector $v$.
    \item $v^Tl \neq 0$.
    \item $v^T \nabla^3 f(x^*)(v,v) \neq 0$.
\end{enumerate}
Then there is a smooth critical curve in $\mathbb{R}^n \times \mathbb{R}$ passing 
through $(x^*,0)$ tangent to the hyperplane $\mathbb{R}^n \times \{0\}$ with no critical point 
on one side of the hyperplane and two critical points on the other side for each $\mu$. The two critical points are hyperbolic and have stable manifolds 
of dimensions $n-1$ and $n$ respectively.
\end{theorem}
\textit{Proof:} See Theorem 3.4.1 in~\cite{guckenheimer2013nonlinear}. $\hfill\square$

Theorem~\ref{t.3.4.1} considers a simple case: a non-strict saddle point $x^*$ of $f$ whose Hessian has a single zero eigenvalue and satisfies a mild third-order condition. It states that if the choice $l = u \in \mathbb{R}^n$ induces a saddle-node bifurcation at $x^*$, then the choice $l = -u$ will instead eliminate the critical point~$x^*$. We now combine Theorem~\ref{t.3.4.1} with a concept that appears trivial at first: for some point $x_0$, the choice $l = -\nabla f(x_0)$ will create a critical point of $f_l$ at $x_0$.
From Theorem~\ref{t.isolatedbifurcations}, this critical point at $x_0$ can only be the result of a bifurcation that occurred in $\Theta(x^*)$, which contains only $x^*$ as a critical point. From Theorem~\ref{t.3.4.1}, if the choice of $l = \nabla f(x_0)$ induces a bifurcation of $x^*$, then the choice of $l = -\nabla f(x_0)$ instead destroys the non-strict critical point.

Theorem~\ref{t.3.4.1} and the above discussion
imply that the choice $l = \nabla f(x_0)$ may be a good candidate for our regularization selection rule.
%This proposed rule is presented in Algorithm~\ref{alg1}, where the update law is $g(x)$ when the gradient is large. 
Under this rule, when $g^k(x)$ enters the small-gradient region $L_{\theta}$ at some point $x_0$, $l$ is set to $\nabla f(x_0)$ and the update law is switched to $g_l(x) = x-\gamma (\nabla f(x) + l)$ until $g^k_l (x)$ leaves $L_{\theta}$. Note that because linear regularization does not affect the Hessian, and by extension the Lipschitz constant $L$, $\gamma$ remains unchanged between $g$ and $g_l$. While this method may bear some superficial similarity to ``momentum methods'' such as in~\cite{jin2018accelerated}, this method differs in that (i) $l$ is not time-varying while in $\Theta(x^*)$, and (ii) momentum methods rely on the SSP.

We note that Theorem~\ref{t.3.4.1} 
provides intuition behind this choice of regularization, but does not provide general theoretical guarantees. To do so we next
determine the general cases for which locally linearly regularized gradient descent avoids non-strict saddles. 

%\begin{figure}[htp]
%    \centering
%    \includegraphics[width=6cm]{BifurcationSetSketch.jpeg}
%    \caption{For a critical point $p_l$ of the regularized function, it lies on a smooth curve of equilibria that either 1) extends back to $p$ with $p_l$ resulting as a bifurcation of a degnerate $p$, 2) 1) extends back to $p$ with $p_l$ resulting as a perturbation of a non-degenerate $p$, or 3) resulted as a hidden bifurcation. However, in all cases the line of equilibria remains inside $\Theta(p)$}
%    \label{fig:bifurcation}
%\end{figure}

\section{Exit Condition of $\Theta(x^*)$}
\label{exitcondition}

By construction, a point $x^*_l \in \Theta(x^*)$ is a critical point of $f_l$ if and only if $\nabla f(x^*_l) = -l$. Because a linear regularizer does not affect the Hessian, $\nabla^2 f_l (x^*_l) = \nabla^2 f (x^*_l)$. That is, if $x^*_l$ is a critical point of $f_l$, its convergence behavior is determined by the Hessian of $f$ at $x^*_l$. In order to analyze this, let us stratify $\Theta(x^*)$ based on the properties of its Hessian:

\begin{definition}
For a $C^2$ function $f: \mathbb{R}^n \rightarrow \mathbb{R}$:
\begin{itemize}
    \item $\Lambda^+ = \{x \in \mathbb{R}^n: \lambda_{min}\nabla^2 f(x) > 0\}$
    \item $\Lambda^0 \,= \{x \in \mathbb{R}^n: \lambda_{min}\nabla^2 f(x) = 0\}$
    \item $\Lambda^- = \{x \in \mathbb{R}^n: \lambda_{min}\nabla^2 f(x) < 0\}$.
\end{itemize}
\end{definition}

Note that $\mathbb{R}^n = \Lambda^+ \cup \Lambda^0 \cup \Lambda^-$. 
%According to Lemma~\ref{l.milnor}, $\nabla^2 f(x^*_l)$ is non-singular, implying $x^*_l \notin \Lambda^0$. 
If $x^*_l \in \Lambda^-$, then it is a strict saddle, and $g^k(x)$ will not converge to $x^*_l$, as shown by the following lemma:

\begin{lemma} \label{l.escapesaddles}
Let $x_0 \in \Theta(x^*)$ for some $x^* \in X^*$ and let $l = \nabla f(x_0)$. Let $Y^*_l = \Theta(x^*) \cap \Lambda^- \cap X^*_l$, where $X^*_l$ is the critical set of $f_l$. Let $\epsilon$ be drawn uniformly from the $n$-Ball with radius $\frac{\theta - \|l\|_2}{L}$. Then
\begin{equation}
    \Pr\left(\lim_{k \rightarrow \infty}g_l^k(x_0 + \epsilon) \in Y^*_l\right) = 0.
\end{equation}
\end{lemma}
\textit{Proof:} All elements of $Y^*_l$ are strict saddle points of the function $f_l$. The map $g_l(x)$ is equivalent to gradient descent on $f_l$. Using this information, Corollary 9 in~\cite{lee2016gradient} provides the result. $\hfill\square$

Note that the one-time perturbation of $x_0$ is done to satisfy a genericity condition necessary to use Corollary~9 in~\cite{lee2016gradient}, and this perturbation is only done when entering $L_\theta$, see~\cite{jin2017escape}. The restriction $\|\epsilon\|_2 \leq \frac{\theta - \|l\|_2}{L}$ ensures $x_0 + \epsilon \in \Theta(x^*)$. Locally linearly regularized gradient descent with this perturbation is presented in Algorithm~1.

\begin{algorithm}[tb]
   \caption{Locally Linearly Regularized Gradient Descent}
   \label{alg1}
\begin{algorithmic}
   \STATE {\textbf{Input}: Stepsize~$\gamma > 0$, Small gradient parameter~$\theta > 0$}
   \FOR{k = 0,1...}
   \IF{$\| \nabla f(x_k)\| > \theta$}
   \STATE {$x_{k+1} \gets x_k - \gamma \nabla f(x_k)$}
   \ELSIF {$\| \nabla f(x_k) \| \leq \theta \And \| \nabla f(x_{k-1})\| > \theta$}
   \STATE {$l \gets \nabla f(x_k)$}
   \STATE {$x_k \gets x_k + \epsilon$ \hfill $\epsilon$ uniformly $\sim \mathbb{B}_0 (\frac{\theta - \|l\|_2}{L})$}
   \STATE {$x_{k+1} \gets x_k - \gamma (\nabla f(x_k) + l)$}
   \ELSE 
   \STATE {$x_{k+1} \gets x_k - \gamma (\nabla f(x_k) + l)$}
   \ENDIF
   \ENDFOR
\end{algorithmic}
\end{algorithm}

%\red{
%Note that, strictly speaking, the initial condition~$x_0$ must be ``generic'' with respect to the gradient field $\nabla f + l$ in order to leverage Corollary~9 from~\cite{lee2016gradient}, which may not be the case with $x_0$ since $l = \nabla f(x_0)$. However, this is easily remedied if, after choosing $l = \nabla f(x_0)$, we perturb to some new point $x_0 + \epsilon$, where $\epsilon$ is a random noise term with $\epsilon \leq \frac{\theta - \|l\|_2}{L}$ to ensure it remains in $\Theta(x^*)$. See~\cite{jin2017escape}.
%}
%\mhmargin{This reads somewhat like a confession...do we need to have noise in our results? If so, I'm good with just doing that. Let's talk about this one.}

If $x^*_l \notin \Lambda^-$, then it must lie in either $\Lambda^0$ or $\Lambda^+$. If $x^*_l \in \Lambda^0$, then $f_l$ does not satisfy the SSP. If $x^*_l \in \Lambda^+$ and $x^*$ is a saddle point, then $x^*_l$ is a \textit{false minimum} by Definition~\ref{d.falsemin}.
%Because $x^*_l \notin \Lambda^0$ and $g_l(x)$ does not converge to any $x^*_l \in \Lambda^- \cap \Theta(x^*)$, then the only way $g_l(x)$ can fail to escape $\Theta(x^*)$ is if $x^*_l \in \Lambda^+ \cap \Theta(x^*)$, which is precisely when $x^*_l$ is a local minimum. Recall that if $x^*$ is not a local minimum of $f$, then such a point~$x^*_l$ is a \textit{false minimum}. 
Therefore, in order to guarantee Algorithm~1 escapes $\Theta(x^*)$ when $x^*$ is a saddle point, we wish to 
show that the choice $l = \nabla f(x_0)$ for $x_0 \in \Theta(x^*)$ always results in $x^*_l \in \Lambda^-$, if $x^*_l$ exists. We formalize this notion with the following definition and assumption:

\begin{definition}
Let~$\Psi(\Theta(x^*)) = \{x \in \Theta(x^*) : \exists y \in \Theta(x^*) \textit{ such that } \nabla f(y) = -\nabla f(x) \textit{ and } y \notin \Lambda^-\}$.
\end{definition}

\begin{assumption} \label{a.nomin}
For the function $f$, for any saddle point~$x^*$,
$\Psi(\Theta(x^*)) = \emptyset$.
\end{assumption}

If $\Psi(\Theta(x^*))$ is nonempty and $x_0 \in \Psi(\Theta(x^*))$, then the choice $l = \nabla f(x_0)$ creates a false minimum or degenerate point in $\Theta(x^*)$. Assumption~\ref{a.nomin} therefore implies that for any saddle point $x^*$ of $f$ and for any point $x_0 \in \Theta(x^*)$, the choice $l = \nabla f(x_0)$ will not create a false minimum or degenerate point in $\Theta(x^*)$. This leads to the main theorem of this work, which addresses the ability of linearly regularized
gradient descent to exit the small-gradient neighborhood of non-strict saddle points in finite time: 

\begin{theorem} \label{t.escapetheta}
Let $x^* \in X^*$ be a saddle point of $f$, and let Assumptions~\ref{a.thetaseparate} and~\ref{a.nomin} hold. 
Let $l = \nabla f(x_0)$ for some $x_0 \in \Theta(x^*)$ with $\theta < \bar{\theta}$. 
Then there almost always exists a finite integer $k_p$ such that $g^{k_{p}}_l(x_0 \red{+ \epsilon}) \notin \Theta(x^*)$.
\end{theorem}
\textit{Proof:} See Appendix~\ref{apx.t.escapetheta}. $\hfill\square$

Theorem~\ref{t.escapetheta} states that under Assumption~\ref{a.nomin}, Algorithm~1 exits $\Theta(x^*)$ for any saddle point $x^*$ in finite time. Note Assumption~\ref{a.nomin} only applies to saddle points, as we do not wish to escape $\Theta(x^*)$ if $x^*$ is a local minimum of $f$. 

Assumption~\ref{a.nomin} gives a sufficient condition for which this regularization method avoids saddles. It is weaker than the SSP, allowing for a class of non-strict saddles. Identifying functions that satisfy Assumption~\ref{a.nomin} is therefore no harder than identifying those with the SSP, and in the following corollaries we identify two properties non-strict saddles may have that are sufficient to satisfy Assumption~\ref{a.nomin}.

\begin{corollary} \label{c.halfspace}
Let~$\nabla f(\Theta(x^*))$ denote the set of all gradients that exist on $\Theta(x^*)$. If $\nabla f(\Theta(x^*))$ lies on an open half-space of $\mathbb{R}^n$, then $\Psi(\Theta(x^*)) = \emptyset$.
\end{corollary}

Trivially, if $x_0 \in \Theta(x^*)$, then $\nabla f(x_0) \in \nabla f(\Theta(x^*))$. If $f$ satisfies the 
condition in Corollary~\ref{c.halfspace}, then $-\nabla f(x_0) \notin \nabla f(\Theta(x^*))$. 
That is, for~$l = \nabla f(x_0)$ no point $x^*_l \in \Theta(x^*)$ exists such that $\nabla f(x^*_l) = -l$, 
which implies $f_l$ has no critical points in $\Theta(x^*)$. Clearly, if $f_l$ has no critical points 
in $\Theta(x^*)$, then
Algorithm~1 exits~$\Theta(x^*)$ by Theorem~\ref{t.escapetheta}. 
Heuristically, if a function can be approximated by an odd polynomial along at least one direction in $\Theta(x^*)$, then by Corollary~\ref{c.halfspace} typically $\Psi(\Theta(x^*)) = \emptyset$, as in Example~\ref{e.halfspace}.
%\mhmargin{Do we have a proof of this, or is this an intuitive thing that seems true?}
%[There is no proof here, this is just me saying whenever I found an example that satisfied this, it was an isolated point with a direction that ``looks'' like $x^3$]

\begin{figure}[htp]
    \centering
    \includegraphics[width=8cm]{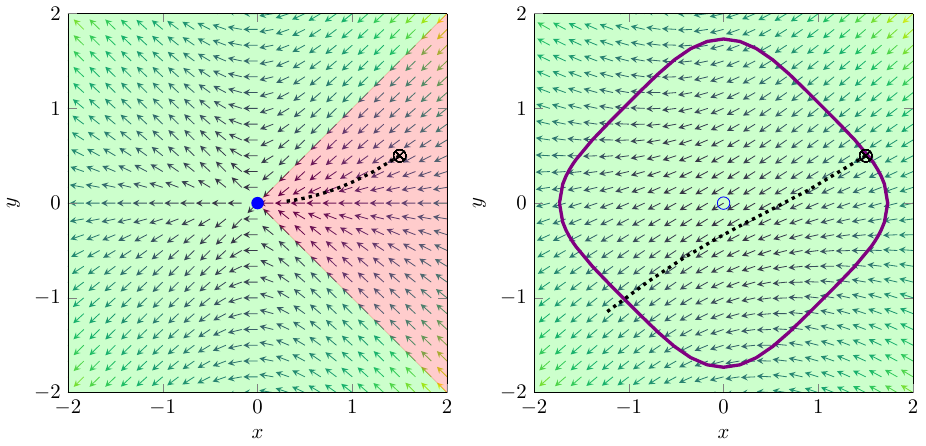}
    \caption{The point $x_0 = (1.5,0.5)$ lies in $W^s_{g}(0,0)$ for the function $f(x,y) = \frac{1}{3}x^3+xy^2$, so $g^k(x_0)$ converges to $(0,0)$ (left). With $l = \nabla f(x_0)$ the critical point at $(0,0)$ is eliminated, and $g_l^k(x_0)$ escapes $\Theta(x_0)$ for $\theta = 3$ in 7 iterations, and enters $W^u_{g}(0,0)$ (right).}
    \label{fig:halfspace}
\end{figure}

\begin{example} \label{e.halfspace}
Consider the function $f(x,y) = \frac{1}{3}x^3+xy^2$, which has a non-strict saddle at $(0,0)$ that satisfies the condition in 
Corollary~\ref{c.halfspace}. This is because $\nabla_x f(x,y) = x^2+y^2$, which is non-negative everywhere. $W^s_g (0,0)$ is represented by the
red region in Figure~\ref{fig:halfspace}, and for every~$x_0 \in W^s(0,0)$, 
we see that the regularzer $l = \nabla f(x_0)$ 
results in no critical points of $f_l$ in $\Theta(0,0)$,
and Algorithm~1 exits $\Theta(0,0)$.
\end{example}

\begin{corollary} \label{c.onlysaddles}
If $\Lambda^- \cap \Theta(p) = \Theta(p)$ then $\Psi(\Theta(p)) = \emptyset$.
\end{corollary}

From Theorem~\ref{t.escapetheta}, if there are only strict saddles in $\Theta(x^*)$ after regularization, then Algorithm~1 exits~$\Theta(x^*)$. 
%False minima must lie in $\Lambda^+$ and if $\Lambda^+ \cap \Theta(x^*)$ is empty (i.e., the Hessian is not positive definite anywhere on $\Theta(x^*)$), then no false minima can exist for any choice of $l$. Any 
Under Corollary~\ref{c.onlysaddles}, critical points of $f_l$ must be strict saddles. Generally, this condition is satisfied by objectives with non-isolated non-strict saddle points, such as in Example~\ref{e.nonpositive}.

\begin{figure}[htp]
    \centering
    \includegraphics[width=8cm]{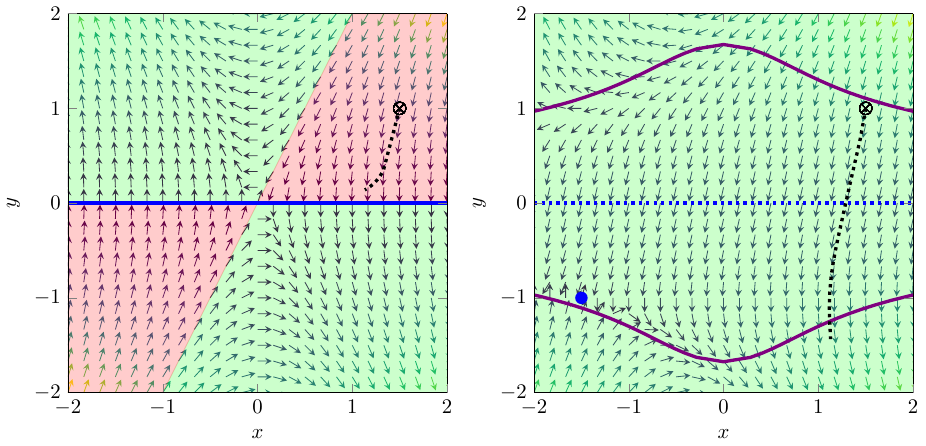}
    \caption{ The function $f(x,y) = \frac{1}{3}xy^3$ has a critical subset on the line $y = 0$. The point $x_0 = (1.5,1)$ lies in $W^s_{g}$ (where $y=0$), so $g^k(x_0)$ converges to $y=0$ (left). With $l = \nabla f(x_0)$, the critical subset at $y = 0$ is eliminated and a strict saddle point of $f_l$ is created at $(-1.5,-1)$. Then Algorithm~1 exits $\Theta(x_0)$ for $\theta = 4.7$ in 15 iterations, and enters $W^u_{g}$ (where~$y=0$) (right).}
    \label{fig:nonpositive}
\end{figure}

\begin{example} \label{e.nonpositive}
Consider the function $f(x,y) = \frac{1}{3}xy^3$, which has a non-strict critical subspace on the $x$-axis. For this function $-\nabla f(x,y) = \nabla f(-x,-y)$, meaning choosing $l = \nabla f(x_0,y_0)$ for \textit{any} $(x_0,y_0)$ will create a critical point of $f_l$ at $(-x_0,-y_0)$. However, $\lambda_{min}(\nabla^2 f (x,y)) < 0$ everywhere with $y \neq 0$, meaning $(-x_0,-y_0)$ will be a strict saddle, and Algorithm~1 exits $\Theta(0,0)$, shown in Figure~\ref{fig:nonpositive}.
\end{example}

\section{The Role of the Hyperparameter $\theta$} \label{hyperparameter}
The behavior of a locally linearly regularized algorithm is highly dependent on the hyperparameter $\theta$. Due to space constraints, determining the upper bound $\bar{\theta}$ from Assumption~\ref{a.thetaseparate} for a particular function $f$ is deferred to a future publication. However, we do wish to illustrate the performance tradeoff between speed and accuracy governed by the choice of $\theta$. Intuitively, small values of $\theta$ should lead to small regularization error. This is formalized in the following theorem.
\begin{theorem} \label{t.error}
Assume $\theta$ is chosen small enough such that, for every critical point~$x^*$ of $f$ that satisfies
$x^* \in \Lambda^+$, we also have $\Theta(x^*) \subset \Lambda^+$. If $\|l\|_2 < \theta$, then $f_l$ will have exactly one critical point $x^*_l$ in $\Theta(x^*)$, and $x^*_l$ will be a non-degenerate minimum. Additionally, if $f$ is $\alpha$-strongly convex on $\Theta(x^*)$, then 
the cost error between~$x^*_l$ and~$x^*$ induced by
regularizing is bounded by
$f(x^*_l) - f(x^*) \leq \frac{\theta^2}{2 \alpha}$.
\end{theorem}

\textit{Proof:} See Appendix~\ref{apx.t.error}. $\hfill\square$

\begin{figure}[htp]
    \centering
    \includegraphics[width=8cm]{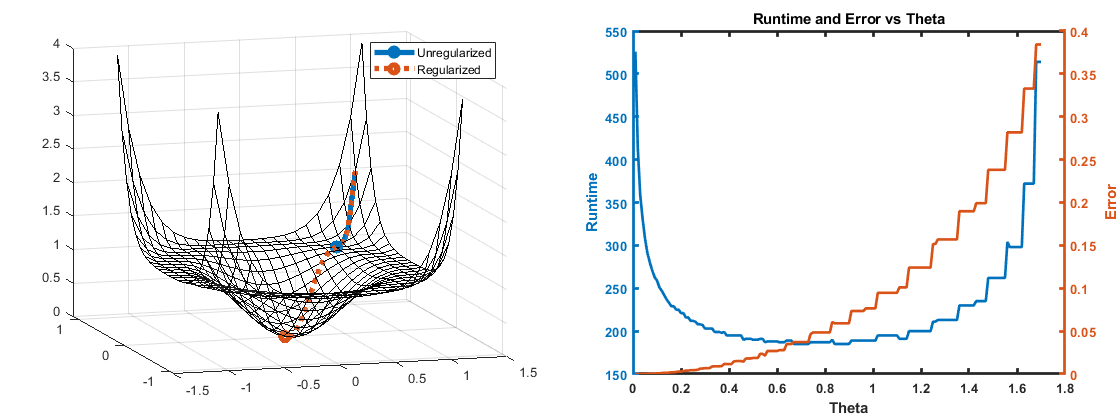}
    \caption{\textit{Left:} Unregularized gradient descent (blue line) converges to the non-strict saddle ring of the inverted wine bottle $f(x,y) = 1+(\sqrt{x^2+y^2}-1)^3 (\sqrt{x^2+y^2}+1)^3$. Algorithm~1 with $\theta = 0.7$ (orange dashed line) converges with minor error. \textit{Right:} Runtime (blue) and final cost error (orange) as $\theta$ is varied. Unregularized gradient descent (not shown) had a final cost error of~$1$ and a runtime of $10,979$.}
    \label{fig:onlyone}
\end{figure}

%\mhmargin{There's a jump in the narrative flow that's kind of jarring.}
The assumption that $f$ is $\alpha$-strongly convex in the neighborhood of local minima is standard in the SSP literature, see Assumption A3.a in~\cite{jin2017escape}. To examine the tradeoff between this error and runtime, we examine the Inverted Wine Bottle, the two-dimensional version of the function in Example~\ref{t.timevaryingbad}. This function has a global minimum at $(0,0)$ surrounded by a ring of non-strict saddles on the unit circle. Unregularized gradient descent initialized outside the unit circle will become stuck and fail to reach the minimum, but locally linearly regularized gradient descent will bypass the ring and reach the origin within some regularization error. We initialize Algorithm~1 at $(1,1)$ with $\gamma = \frac{1}{54}$ and run using values of $\theta$ varying from $0.01$ to $1.7$ ($\bar{\theta} \approx 1.717$ for this function). Each run of the algorithm terminates when $\|\nabla f(x) + l\| \leq 10^{-7}$. The runtime and final cost error due to regularization are plotted in Figure~\ref{fig:onlyone}.

Figure~\ref{fig:onlyone} shows that final cost error increases with $\theta$, as expected from Theorem~\ref{t.error}, but the relationship between $\theta$ and the runtime is more complex. Initially, as $\theta$ is varied away from $0$, the runtime decreases. This is intuitive, as smaller choices of $\theta$ limit the use of regularizers to smaller regions of the  space of iterates. However, as $\theta$ approaches $\bar{\theta}$, the runtime increases. This is due to the large perturbation of the minimum resulting from the large value of $l$. That is, for small values of $\theta$ the algorithm takes a long time to escape saddle points, and for large values of $\theta$ it takes a long time to converge to the minimum. A full analysis of how to tune $\theta$ and its effects on the performance of a locally linearly regularized algorithm is the subject of future work.

\section{Concluding Remarks}
\label{concludingremarks}

We have answered Question~\ref{q} by demonstrating that linear regularizers can be used to enforce the SSP for non-convex objective functions, and that any such regularization scheme must  both do so locally and must choose $l$ based on first-order information. We have presented a local linear regularization scheme with these properties that enforces satisfaction of the SSP. This scheme is proven to escape a broad class of isolated and non-isolated non-strict saddle points. Future work will address tuning the hyperparameter $\theta$.
%, and \red{establishing a proof of convergence}}.
%\mhmargin{What if we remove the red part here? It's too easy for a reader to get here, see this, and think ``they haven't done anything!'' We have, but we don't yet have
%a proof of convergence to minimizers (which is fine). I just don't want a reader to get to the end and think we've just been messing around rather than showing
%anything serious. Let me know what you think.}

\appendix
\section{Appendix}
\subsection{Proof of Theorem~\ref{t.isolatedbifurcations}} \label{apx.t.isolatedbifurcations}

Consider the function $h(x,\mu) = \nabla f(x) + \mu l$ with $\|l\|_2 < \theta$. $h$ maps $\mathbb{R}^n \times \mathbb{R} \rightarrow \mathbb{R}^n$, and $(x^*,0)$ represents a critical point of the non-regularized function $f$. Consider a point $(x^*_l,1) \in \mathbb{R}^n \times [0,1]$ where $\nabla f(x^*_l) + l = 0$, which corresponds to a critical point of the regularized function $f(x) + l^T x$. 
%From Lemma~\ref{l.milnor}, $\nabla^2 f(x^*_l)$ is non-singular almost always. 
If the critical point of $f_l$ at $x^*_l$ resulted as a bifurcation originating at $x^*$, then the Implicit Function Theorem (Theorem 2.3 in~\cite{matsumoto2002introduction}) states that there exists an open neighborhood $U \subset \mathbb{R}$ containing $\mu = 1$ such that there exists a smooth function $\zeta: U \rightarrow \mathbb{R}^n$ such that $\zeta(1) = x^*_l$ and $\nabla f(\zeta(\mu)) + \mu l = 0$ for all $\mu \in U$. That is, starting at $\mu = 1$ and moving in the negative direction, $(\zeta(\mu), \mu)$ is a smooth curve in $\mathbb{R}^n \times [0,1]$ that describes the location of a critical point for different values of $\mu$. 
%Observe that $x^*_l$ is due to either a bifurcation of~$x^*$ (if $x^*$ is degenerate) or a perturbation of~$x^*$ (if $x^*$ is non-degenerate) if $(x^*,0)$ and $(x^*_l,1)$ are connected by the smooth curve defined by $\zeta$. 
Because $\|\nabla f(\zeta(\mu))\|_2 = \mu\|l\|_2 < \theta$ for all $\mu \in [0,1]$, this curve must lie in the connected subset of $L_{\theta} \times [0,1]$ that contains $(x^*,0)$, which is $\Theta(x^*) \times [0,1]$. Therefore $x^*_l \in \Theta(x^*)$. $\hfill\square$

\subsection{Proof of Theorem~\ref{t.escapetheta}} \label{apx.t.escapetheta}

The map $g_l(x)$ is equivalent to gradient descent on the function $f_l(x) = f(x) + l^Tx$. Under Assumption~\ref{a.nomin}, any critical points in $\Theta(x^*)$ must lie in $\Lambda^-$, which implies they are strict saddles. Lemma~\ref{l.escapesaddles} states $\lim_{k \rightarrow \infty}g_l^k(x_0)$ is almost never a strict saddle. Therefore $\lim_{k \rightarrow \infty}g_l^k(x_0)$ will almost always lie outside $\Theta(x^*)$, implying it exits $\Theta(x^*)$ in finite time. $\hfill\square$

\subsection{Proof of Theorem~\ref{t.error}} \label{apx.t.error}
From Theorem~\ref{t.isolatedbifurcations}, a perturbation of $x^*$ remains in $\Theta(x^*)$. Because $\Lambda^0 \cap \Theta(x^*) = \emptyset$, no point $ x \in \Theta(x^*)$ has $\nabla^2f(x)$ singular, therefore there exists exactly one point $x^*_l \in \Theta(x^*)$ for which $\nabla f(x^*_l) + l= 0$, and $x^*_l \in \Lambda^+$. $\alpha$-strong convexity on $\Theta(x^*)$ implies that, for every point $x \in \Theta(x^*)$, $\frac{1}{2}\|\nabla f(x)\|^2 \geq \alpha (f(x) - f(x^*))$
holds. Since $x^*_l \in \Theta(x^*)$, then $\| \nabla f(x^*_l) \|_2 \leq \theta$. The result follows by substitution. $\hfill\square$

% Use \bibliography{yourbibfile} instead or the References section will not appear in your paper
\bibliography{Biblio}

%\section{Acknowledgments}
%AAAI is especially grateful to Peter Patel Schneider for his work in implementing the original aaai.sty file, liberally using the ideas of other style hackers, including Barbara Beeton. We also acknowledge with thanks the work of George Ferguson for his guide to using the style and BibTeX files --- which has been incorporated into this document --- and Hans Guesgen, who provided several timely modifications, as well as the many others who have, from time to time, sent in suggestions on improvements to the AAAI style. We are especially grateful to Francisco Cruz, Marc Pujol-Gonzalez, and Mico Loretan for the improvements to the Bib\TeX{} and \LaTeX{} files made in 2020.

%The preparation of the \LaTeX{} and Bib\TeX{} files that implement these instructions was supported by Schlumberger Palo Alto Research, AT\&T Bell Laboratories, Morgan Kaufmann Publishers, The Live Oak Press, LLC, and AAAI Press. Bibliography style changes were added by Sunil Issar. \verb+\+pubnote was added by J. Scott Penberthy. George Ferguson added support for printing the AAAI copyright slug. Additional changes to aaai23.sty and aaai23.bst have been made by Francisco Cruz, Marc Pujol-Gonzalez, and Mico Loretan.

%\bigskip
%\noindent Thank you for reading these instructions carefully. We look forward to receiving your electronic files!

\end{document}